\documentclass[12pt,oneside,reqno]{amsart}
\usepackage{txfonts}
\usepackage{bbm}
\usepackage{amsmath}
\usepackage{graphicx}
\usepackage{mathrsfs}
\usepackage{stmaryrd}
\usepackage{amsfonts}
\usepackage{enumerate,amsmath,amssymb,amsthm}

\numberwithin{equation}{section}

\newcommand{\be}{\begin{eqnarray}}
\newcommand{\ee}{\end{eqnarray}}
\newcommand{\ce}{\begin{eqnarray*}}
\newcommand{\de}{\end{eqnarray*}}
\newtheorem{theorem}{Theorem}[section]
\newtheorem{lemma}[theorem]{Lemma}
\newtheorem{remark}[theorem]{Remark}
\newtheorem{definition}[theorem]{Definition}
\newtheorem{proposition}[theorem]{Proposition}
\newtheorem{Examples}[theorem]{Example}
\newtheorem{corollary}[theorem]{Corollary}

\def\v{{\mathbf{v}}}
\def\eps{\varepsilon}

\def\p{\partial}

\def\[{{\Big[}}
\def\]{{\Big]}}
\def\<{{\langle}}
\def\>{{\rangle}}
\def\({{\Big(}}
\def\){{\Big)}}

\def\bx{{\mathbf{x}}}

\def\dif{{\mathord{{\rm d}}}}

\def\no{\nonumber}
\def\={&\!\!=\!\!&}
\def\bt{\begin{theorem}}
\def\et{\end{theorem}}
\def\bl{\begin{lemma}}
\def\el{\end{lemma}}
\def\br{\begin{remark}}
\def\er{\end{remark}}

\def\bd{\begin{definition}}
\def\ed{\end{definition}}
\def\bp{\begin{proposition}}
\def\ep{\end{proposition}}
\def\bc{\begin{corollary}}
\def\ec{\end{corollary}}
\def\bx{\begin{Examples}}
\def\ex{\end{Examples}}

\def\cF{{\mathcal F}}

\def\cL{{\mathcal L}}
\def\cM{{\mathcal M}}

\def\cO{{\mathcal O}}

\def\cS{{\mathcal S}}
\def\cT{{\mathcal T}}

\def\mD{{\mathbb D}}
\def\mE{{\mathbb E}}

\def\mH{{\mathbb H}}

\def\mN{{\mathbb N}}

\def\mR{{\mathbb R}}
\def\mS{{\mathbb S}}

\def\mW{{\mathbb W}}

\def\sB{{\mathscr B}}

\def\sF{{\mathscr F}}

\def\sS{{\mathscr S}}
\def\sT{{\mathscr T}}

\def\geq{\geqslant}
\def\leq{\leqslant}

\allowdisplaybreaks

\begin{document}

\title{Discontinuous Stochastic Differential Equations Driven by L\'evy Processes}

\date{}
\author{Xicheng Zhang}


\dedicatory{
School of Mathematics and Statistics,
Wuhan University, Wuhan, Hubei 430072, P.R.China,\\
Email: XichengZhang@gmail.com
 }

\begin{abstract}
In this article we prove the pathwise uniqueness for stochastic differential equations in $\mR^d$
with time-dependent Sobolev drifts, and driven by symmetric $\alpha$-stable processes
provided that $\alpha\in(1,2)$ and its spectral measure is non-degenerate.
In particular, the drift is allowed to have jump discontinuity when
$\alpha\in(\frac{2d}{d+1},2)$. Our proof is based on some estimates of Krylov's type
for purely discontinuous semimartingales.
\end{abstract}

\maketitle
\rm

\section{Introduction and Main Result}

Consider the following SDE driven by a symmetric $\alpha$-stable  process in $\mR^d$:
\begin{align}
\dif X_t=b_t(X_t)\dif t+\dif L_t,\ \ X_0=x,\label{SDE}
\end{align}
where $b:\mR_+\times\mR^d\to\mR^d$ is a measurable function, $(L_t)_{t\geq 0}$ is a $d$-dimensional symmetric
$\alpha$-stable process defined on some filtered probability space $(\Omega,\sF,P;(\sF_t)_{t\geq 0})$.
The aim of this paper is to study the pathwise uniqueness of SDE (\ref{SDE}) with discontinuous $b$.

Let us first briefly recall some well known results in this direction.
When $L_t$ is a standard $d$-dimensional Brownian motion,
Veretennikov \cite{Ve} first proved the existence of a unique strong solution for SDE (\ref{SDE})
with bounded measurable $b$. In \cite{Kr-Ro}, Krylov and R\"ockner relaxed the boundedness assumptions on
$b$ to the following integrability assumptions:
\begin{align}
\int^T_0\left(\int_{\mR^d}|b_t(x)|^p\dif x\right)^{\frac{q}{p}}\dif t<+\infty, \ \ \forall T>0,\label{BP3}
\end{align}
provided that
\begin{align}
\frac{d}{p}+\frac{2}{q}<1.\label{BP4}
\end{align}
Recently, in \cite{Zh0} we extended Krylov and R\"ockner's result to the case of non-constant Soboev diffusion coefficients
and meanwhile, obtained the stochastic homeomorphism flow property of solutions and the strong Feller property.

In the case of symmetric $\alpha$-stable processes, the pathwise uniqueness for SDE (\ref{SDE}) with irregular drift
is far from being complete. When $d=1, \alpha\in[1,2)$ and
$b$ is time-independent and bounded continuous, Tanaka, Tsuchiya and Watanabe \cite{Ta-Ts-Wa} proved
the pathwise uniqueness of solutions to SDE (\ref{SDE}). When $d>1, \alpha\in[1,2)$,
the spectral measure of $L_t$ is non-degenerate,
and $b$ is time-independent and bounded H\"older continuous, where the H\"older index $\beta$ satisfies
$$
\beta>1-\frac{\alpha}{2},
$$
Priola \cite{Pr} recently proved the pathwise uniqueness and the stochastic homeomorphism flow property of solutions to SDE (\ref{SDE}).
When $d=1,\alpha\in(1,2)$ and $b$ is only bounded measurable, Kurenok \cite{Ku} obtained the existence of
weak solutions for SDE (\ref{SDE}) by proving an estimate of Krylov's type: for any $T>0$,
\begin{align}
\mE\int_0^Tf_t(X_t)\dif t\leq C\|f\|_{L^2([0,T]\times\mR)}.\label{Kry}
\end{align}
When $d>1$, $\alpha\in(1,2)$ and $b$ is time-independent and belongs to some Kato's class,
Chen, Kim and Song \cite[Theorem 2.5]{Ch-Ki-So} proved the existence of martingale solutions
(equivalently weak solutions) in terms of Feller semigroup (cf. \cite{Bo-Ja}). On the other hand, there are many
works devoted to the study of weak uniqueness (i.e., the well-posedness of martingale problems) for SDEs with jumps.
This is refereed to the survey paper of Bass \cite{Ba3}.
However, to the author's knowledge, there are few results about the pathwise uniqueness
for multidimensional SDE (\ref{SDE}) with discontinuous drifts.

Before stating our main result, we recall some facts about symmetric $\alpha$-stable processes.
Let $(L_t)_{t\geq 0}$ be a $d$-dimensional symmetric $\alpha$-stable process.
By L\'evy-Khinchin's formula, its characteristic function is given by (cf. \cite{Sa})
$$
\mE e^{i\xi L_t}=e^{-t\psi(\xi)},
$$
where
$$
\psi(\xi)=\int_{\mR^d}(1-e^{i\<\xi,x\>}+i\<\xi,x\>1_{|x|\leq 1})\nu(\dif x),
$$
and the L\'evy measure $\nu$ with $\nu(\{0\})=0$ is given by
\begin{align}
\nu(U)=\int_{\mS^{d-1}}\!\!\int^{+\infty}_0
\frac{1_U(r\theta)}{r^{d+\alpha}}\dif r\mu(\dif\theta),\ \ U\in\sB(\mR^d),\label{Coa}
\end{align}
where $\mu$ is a symmetric finite measure on the unit sphere $\mS^{d-1}:=\{\theta\in\mR^d: |\theta|=1\}$,
called {\it spectral measure} of stable process $L_t$. By an elementary calculation, we have
$$
\psi(\xi)=\int_{\mR^d}(1-\cos\<\xi,x\>)\nu(\dif x)=c_\alpha\int_{\mS^{d-1}}|\<\xi,\theta\>|^\alpha\mu(\dif \theta).
$$
In particular, if $\mu$ is the uniform distribution on $\mS^{d-1}$,
then $\psi(\xi)=c_\alpha|\xi|^\alpha$, here $c_\alpha$ may be different.
Throughout this paper, we make the following assumption:

{\bf (H$^\alpha$)}:  For some $\alpha\in(0,2)$ and constant $C_\alpha>0$,
\begin{align}
\psi(\xi)\geq C_\alpha|\xi|^\alpha,\ \ \forall\xi\in\mR^d.
\end{align}
We remark that the above condition is equivalent that
the support of  spectral measure $\mu$ is not contained in a proper linear subspace of $\mR^d$
(cf. \cite[page 4]{Pr}).

We now introduce the class of local strong solutions for SDE (\ref{SDE}). Let $\tau$ be any ($\sF_t$)-stopping time.
For $x\in\mR^d$, let $\sS^\tau_{b}(x)$ be the class of all $\mR^d$-valued ($\sF_t$)-adapted c\`adl\`ag
stochastic process $X_t$ on $[0,\tau)$ satisfying
$$
P\left\{\omega: \int^T_0|b_s(X_s(\omega))|\dif s<+\infty,
\forall T\in[0,\tau(\omega))\right\}=1,
$$
and such that
$$
X_t=x+\int^t_0b_s(X_s)\dif s+L_t,\ \ \forall t\in[0,\tau),\ \ a.s.
$$

The main result of the present paper is:
\bt\label{Main}
Assume that {\bf (H$^\alpha$)} holds with  $\alpha\in(1,2)$, and $b:\mR_+\times\mR^d\to\mR^d$
satisfies that for some $\beta\in(1-\frac{\alpha}{2},1)$, $p>\frac{2d}{\alpha}$
and any $T,R>0$,
\begin{align}
\sup_{t\in[0,T]}\int_{B_R}\!\int_{B_R}\frac{|b_t(x)-b_t(y)|^p}{|x-y|^{d+\beta p}}\dif x\dif y<+\infty\label{Op6}
\end{align}
and
\begin{align}
\sup_{(t,x)\in[0,T]\times B_R}|b_t(x)|<+\infty,\label{Op06}
\end{align}
where $B_R=\{x\in\mR^d: |x|\leq R\}$.
Then, for any $x\in\mR^d$, there exists an ($\sF_t$)-stopping time $\zeta(x)$ (called explosion time) and
a unique strong solution $X_t\in\sS^{\zeta(x)}_b(x)$ to SDE (\ref{SDE}) with
\begin{align}
\lim_{t\uparrow\zeta(x)}X_t(x)=+\infty,\ \ a.s.\label{HH1}
\end{align}
\et
\br
Let $\cO$ be a bounded smooth domain in $\mR^d$.
It is well known that for any $\beta\in(0,1)$ and $p\in(1,\frac{1}{\beta})$ (cf. \cite{Ar-Zh}),
\begin{align}
1_\cO\in\mW^\beta_p,\label{Fact}
\end{align}
where $\mW^\beta_p$ is the fractional Sobole space defined by (\ref{LL1}) below.
Hence, if $\alpha\in(\frac{2d}{d+1},2)$, then one can choose
$$
\beta\in(1-\frac{\alpha}{2},\frac{\alpha}{2d}),\ \ p\in(\frac{2d}{\alpha},\frac{1}{\beta})
$$
so that Theorem \ref{Main} can be used to uniquely solve the following discontinuous SDE:
$$
\dif X_t=[b^{(1)}1_\cO+b^{(2)}1_{\cO^c}](X_t)\dif t+\dif L_t,\ \ X_0=x,
$$
where $b^{(i)}, i=1,2$ are two bounded and locally H\"older continuous functions
with H\"older index greater than $\beta$. In one dimensional case, if $\alpha\in(1,2)$, it is
well known that regularity condition (\ref{Op6}) can be dropped (cf. \cite[p.82, Remark 1]{Ta-Ts-Wa}).
The key point in this case is that the weak uniqueness is equivalent to the pathwise uniqueness.
However, in the case of $d\geq 2$,
it is still open that whether SDE (\ref{SDE}) has a unique strong solution when $b$ is only bounded measurable.
\er

For proving this theorem, as in \cite{Zv,Kr-Ro,Zh0}, we mainly study the following partial integro-differential equation
(abbreviated as PIDE) by using some interpolation techniques:
$$
\p_t u=\cL_0u+b^i\p_i u+f, \ \ u_0(x)=0,
$$
where $\cL_0$ is the generator of L\'evy process $(L_t)_{t\geq 0}$ given by
\begin{align}
\cL_0 u(x)=\int_{\mR^d}(u(x+z)-u(x)-1_{|z|\leq 1}z^i\p_i u(x))\nu(\dif z)=
\lim_{\eps\downarrow 0}\int_{|z|>\eps}(u(x+z)-u(x))\nu(\dif z),\label{Op4}
\end{align}
where the second equality is due to the symmetry of $\nu$.
Here and below, we use the convention that the repeated index will be summed automatically.
However, we need to firstly extend Krylov's estimate (\ref{Kry}) to the multidimensional case.
As in \cite{Ku}, we shall investigate the following semi-linear PIDE:
$$
\p_t u=\cL_0u+\kappa|\nabla u|+f, \ \ u_0(x)=0,
$$
where $\kappa>0$ and $\nabla$ is the gradient operator with respect to the spatial variable $x$.
We want to emphasize that Fourier's transform used in \cite{Ku1,Ku} seems only work for one-dimensional case.

Our method for studying the above two PIDEs is based on semigroup arguments. For this aim, we shall derive some smoothing
and asymptotic properties about the Markovian semigroup associated with $\cL_0$. In particular, the interpolation
techniques will be used frequently. This will be done in Section 2.
In Section 3, partly following Kurenok's idea, we shall prove two Krylov's estimates  for
multidimensional purely discontinuous semimartingales. In Section 4, we prove our main result by using Zvonkin's
transformation of phase space to remove the drift.

In the remainder of this paper, the letter $C$ with or without subscripts will denote a positive constant
whose value may change in different occasions.
\section{Preliminaries}

For $p\geq 1$, the norm in $L^p$-space $L^p(\mR^d)$ is denoted by $\|\cdot\|_p$. For $\beta\geq 0$ and $p\geq 1$, let $\mH^{\beta}_p$
be the space of Bessel potential, i.e.,
$$
\mH^{\beta}_p=(I-\Delta)^{-\beta}(L^p(\mR^d)).
$$
In other words, $\mH^{\beta}_p$ is the domain of fractional operator $(I-\Delta)^\beta$, where
$(I-\Delta)^\beta$ is defined through
$$
(I-\Delta)^\beta f =\cF^{-1}((1+|\cdot|^2)^\beta(\cF f)),\ \ f\in C^\infty_0(\mR^d),
$$
where $\cF$ (resp. $\cF^{-1}$) denotes the Fourier transform (resp. the Fourier inverse transform).
Notice that for $\beta=m\in\mN$, an equivalent norm of $\mH^{\beta}_p$ is given by (cf. \cite[p.177]{Tr})
$$
\|f\|_{m,p}=\|f\|_p+\|\nabla^m f\|_p.
$$
By Sobolev's embedding theorem, if $\beta-\frac{d}{p}>0$ is not an integer, then (cf. \cite[p.206, (16)]{Tr})
\begin{align}
\mH^{\beta}_p\hookrightarrow C^{\beta-\frac{d}{p}}(\mR^d),\label{em}
\end{align}
where for $\gamma>0$, $C^\gamma(\mR^d)$ is the usual H\"older space with the norm:
$$
\|f\|_{C^\gamma}:=\sum_{k=0}^{[\gamma]}\sup_{x\in\mR^d}|\nabla^k f(x)|+
\sup_{x\not=y}\frac{|\nabla^{[\gamma]} f(x)-\nabla^{[\gamma]} f(y)|}{|x-y|^{\gamma-[\gamma]}},
$$
where $[\gamma]:=\max\{m\in\mN: m\leq\gamma\}$ is the integer part of $\gamma$.

Let $A$ and $B$ be two Banach spaces. For $\theta\in[0,1]$, we use $[A,B]_{\theta}$
to denote the complex interpolation space between $A$ and $B$.
We have the following relation (cf. \cite[p.185, (11)]{Tr}):
for $p>1$, $\beta_1\not=\beta_2$ and $\theta\in(0,1)$,
\begin{align}
[\mH^{\beta_1}_p,\mH^{\beta_2}_p]_{\theta}=\mH^{\beta_1+\theta(\beta_2-\beta_1)}_p.\label{Sob}
\end{align}
On the other hand, for $0<\beta\not=$integer, the fractional Sobolev space $\mW^\beta_p$ is defined by
(cf. \cite[p.190,(15)]{Tr})
\begin{align}
\|f\|^\sim_{\beta,p}:=\|f\|_p+\sum_{k=0}^{[\beta]}\left(\int_{\mR^d}\!\!
\int_{\mR^d}\frac{|\nabla^k f(x)-\nabla^k f(y)|^p}{|x-y|^{d+(\beta-[\beta])p}}\dif x\dif y\right)^{1/p}<+\infty.\label{LL1}
\end{align}
For $\beta=0,1,2,\cdots$, we set $\mW^\beta_p:=\mH^\beta_p$. The relation between $\mH^\beta_p$ and $\mW^\beta_p$
is given as follows (cf. \cite[p.180, (9)]{Tr}): for any $\beta>0$, $\eps\in(0,\beta)$ and $p>1$,
\begin{align}
\mH^{\beta+\eps}_p\hookrightarrow \mW^\beta_p\hookrightarrow\mH^{\beta-\eps}_p.\label{Fr}
\end{align}

We recall the following complex interpolation theorem (cf. \cite[p.59, Theorem (a)]{Tr}).
\bt\label{Inter}
Let $A_i\subset B_i, i=0,1$ be Banach spaces. Let $\sT:A_i\to B_i, i=0,1$ be bounded linear operators.
For $\theta\in[0,1]$,  we have
$$
\|\sT\|_{A_\theta\to B_\theta}\leq\|\sT\|^{1-\theta}_{A_0\to B_0}\|\sT\|^{\theta}_{A_1\to B_1},
$$
where $A_\theta:=[A_0,A_1]_{\theta}$, $B_\theta:=[B_0,B_1]_{\theta}$, and
$\|\sT\|_{A_\theta\to B_\theta}$ denotes the operator norm of $\sT$ mapping $A_\theta$ to $B_\theta$.
\et
Let $f$ be a locally integrable function on $\mR^d$. The Hardy-Littlewood maximal function is defined by
$$
\cM f(x):=\sup_{0<r<\infty}\frac{1}{|B_r|}\int_{B_r}|f(x+y)|\dif y,
$$
where $B_r:=\{x\in\mR^d: |x|<r\}$.
The following well known results can be found in  \cite{Re-Zh, Zh1} and \cite[page 5, Theorem 1]{St}.
\bl\label{Le2}
(i) For $f\in \mW^1_1$, there exists a constant $C_d>0$ and a Lebesgue zero set $E$ such that for all
$x,y\notin E$,
\begin{align}
|f(x)-f(y)|\leq C_d\cdot |x-y|\cdot(\cM|\nabla f|(x)+\cM|\nabla f|(y)).\label{Es2}
\end{align}
(ii) For $p>1$, there exists a constant $C_{d,p}>0$ such that for all $f\in L^p(\mR^d)$,
\begin{align}
\|\cM f\|_p\leq C_{d,p}\|f\|_p.\label{Es30}
\end{align}
\el
For fixed $z\in\mR^d$, define the shift operator
$$
\sT_zf(x):=f(x+z)-f(x).
$$
We have the following useful estimate.
\bl\label{Le4}
For $p>1$ and $\gamma\in[1,2]$, there exists a constant
$C=C(p,\gamma,d)>0$ such that for any $f\in\mH^\gamma_p$,
\begin{align}
\|\sT_zf\|_{1,p}\leq C|z|^{\gamma-1}\cdot\|f\|_{\gamma,p}.\label{Lp7}
\end{align}
\el
\begin{proof}
By (\ref{Es2}), we have for Lebesgue almost all $x\in\mR^d$,
$$
|\sT_zf(x)|\leq C|z|\cdot(\cM|\nabla f|(x+z)+\cM|\nabla f|(x)),
$$
and so, by (\ref{Es30}),
$$
\|\sT_zf\|_p\leq C|z|\cdot\|\cM|\nabla f|\|_p\leq C|z|\cdot\|\nabla f\|_p\leq C|z|\cdot\|f\|_{1,p}.
$$
On the other hand, it is clear that for any $\beta>0$,
$$
\|\sT_zf\|_{\beta,p}\leq 2\|f\|_{\beta,p}.
$$
By Theorem \ref{Inter} and (\ref{Sob}), for $\theta\in(0,1)$, we immediately have
$$
\|\sT_zf\|_{\theta \beta,p}\leq C|z|^{1-\theta}\cdot\|f\|_{1+\theta(\beta-1),p},
$$
which gives the desired result by letting $\theta=2-\gamma$ and $\beta=\frac{1}{2-\gamma}$.
\end{proof}

We now recall the following well known properties about the symmetric $\alpha$-stable process $(L_t)_{t\geq 0}$
(cf. \cite[Theorem 25.3]{Sa} and \cite[Lemma 3.1]{Pr}).
\bp\label{Pr1}
Let $\mu_t$ be the law of $\alpha$-stable process $L_t$.
\begin{enumerate}[(i)]
\item (Scaling property): For any $\lambda>0$ , $(L_t)_{t\geq 0}$ and $(\lambda^{-\frac{1}{\alpha}}L_{\lambda t})_{t\geq 0}$ have
the same finite dimensional law. In particular, for any $t>0$ and $A\in\sB(\mR^d)$, $\mu_t(A)=\mu_1(t^{-\frac{1}{\alpha}}A)$.
\item (Existence of smooth density): For any $t>0$, $\mu_t$ has a smooth density $p_t$ with
respect to the Lebesgue measure, which is given by
$$
p_t(x)=\frac{1}{(2\pi)^d}\int_{\mR^d}e^{-i\<x,\xi\>}e^{-t\psi(\xi)}\dif\xi.
$$
Moreover, $p_t(x)=p_t(-x)$ and for any $k\in\mN$, $\nabla^kp_t\in L^1(\mR^d)$.
\item (Moments): For any $t>0$, if $\beta<\alpha$, then $\mE |L_t|^\beta<+\infty$;
if $\beta\geq\alpha$, then $\mE |L_t|^\beta=\infty$.
\end{enumerate}
\ep

The Markovian semigroup associated with the L\'evy process $(L_t)_{t\geq 0}$ is given by
\begin{align}
\cT_t f(x)=\mE(f(L_t+x))=\int_{\mR^d}p_t(z-x)f(z)\dif z=\int_{\mR^d}p_t(x-z)f(z)\dif z.\label{Semi}
\end{align}
We have:
\bl\label{L11}
(i) For any $\alpha\in(0,2)$, $p>1$ and $\beta,\gamma\geq 0$, we have for all $f\in \mH^\beta_p$,
\begin{align}
\|\cT_t f\|_{\beta+\gamma,p}\leq C t^{-\gamma/\alpha}\|f\|_{\beta,p}.\label{Lp1}
\end{align}
(ii) For any $\alpha\in(1,2)$, $\theta\in[0,1]$ and $p>1$, there exists a constant $C=C(d,p,\theta)>0$ such that for all
$f\in\mH^\theta_p$,
\begin{align}
\|\cT_t f-f\|_p\leq Ct^{\theta/\alpha}\|f\|_{\theta,p}.\label{Lp2}
\end{align}
\el
\begin{proof}
(i) Let $f\in C^\infty_0(\mR^d)$. For any $k,m\in\mN$, by the scaling property, we have
$$
\nabla^{k+m}\cT_t f(x)=t^{-(d+k)/\alpha}\int_{\mR^d}(\nabla^kp_1)(t^{-1/\alpha}(x-z))\nabla^mf(z)\dif z.
$$
Hence,
$$
\|\nabla^{k+m}\cT_t f\|_p\leq t^{-k/\alpha}\|\nabla^mf\|_p\int_{\mR^d}|\nabla^kp_1|(x)\dif x.
$$
Since $C^\infty_0(\mR^d)$ is dense in $\mH^m_p$, we further have for any $f\in\mH^m_p$,
$$
\|\nabla^{k+m}\cT_t f\|_p\leq t^{-k/\alpha}\|f\|_{m,p}\int_{\mR^d}|\nabla^kp_1|(x)\dif x.
$$
On the other hand, it is clear that
$$
\|\cT_t f\|_p\leq \|f\|_p.
$$
By Theorem \ref{Inter}, we obtain (\ref{Lp1}).

(ii) First, we assume that $f\in \mH^1_p$. By (\ref{Es2}), we have for Lebesgue almost all $x\in\mR^d$,
\begin{align*}
|\cT_tf(x)-f(x)|&\leq \int_{\mR^d}|f(x+y)-f(x)|\cdot p_t(y)\dif y\\
&\leq C\int_{\mR^d}(\cM|\nabla f|(x+y)+\cM|\nabla f|(x))\cdot|y|\cdot p_t(y)\dif y,
\end{align*}
and so, by (\ref{Es30}) and the scaling property,
\begin{align*}
\|\cT_tf-f\|_p&\leq C\|\cM|\nabla f|\|_p\int_{\mR^d}|y|\cdot p_t(y)\dif y\\
&\leq C\|\nabla f\|_p\mE|L_t|=Ct^{1/\alpha}\|\nabla f\|_p\mE|L_1|.
\end{align*}
Estimate (\ref{Lp2}) follows by (iii) of Proposition \ref{Pr1} and Theorem \ref{Inter} again.
\end{proof}
We also need the following simple result for proving the uniqueness.
\bl\label{Le0}
Let $(Z_t)_{t\geq 0}$ be a locally bounded and ($\sF_t$)-adapted process and $(A_t)_{t\geq 0}$
a continuous real valued non-decreasing ($\sF_t$)-adapted process with $A_0=0$.
Assume that for any stopping time $\eta$ and $t\geq 0$,
$$
\mE|Z_{t\wedge\eta}|\leq\mE\int^{t\wedge\eta}_0|Z_s|\dif A_s.
$$
Then $Z_t=0$ a.s. for all $t\geq 0$.
\el
\begin{proof}
By replacing $A_t$ by $A_t+t$, one may assume that $t\mapsto A_t$ is strictly increasing.
For $t\geq 0$, define the stopping time
$$
\tau_t:=\inf\{s\geq 0: A_s\geq t\}.
$$
It is clear that $\tau_t$ is the inverse of $t\mapsto A_t$. Fix $T>0$.
By the assumption and the change of variable, we have
$$
\mE|Z_{T\wedge \tau_t}|\leq\mE\int^{T\wedge \tau_t}_0|Z_s|\dif A_s
\leq\mE\int^{\tau_t}_0|Z_{T\wedge s}|\dif A_s=\int^{t}_0\mE|Z_{T\wedge \tau_s}|\dif s.
$$
By Gronwall's inequality, we obtain $Z_{T\wedge \tau_t}=0$. Letting $t\to\infty$ gives the conclusion.
\end{proof}

\section{Krylov's estimates for purely discontinuous semimartingales}

Let $(L_t)_{t\geq 0}$ be a symmetric $\alpha$-stable process. The associated Poisson random measure is defined by
$$
N((0,t]\times U):=\sum_{s\in(0,t]}1_U(L_s-L_{s-}),\ \ U\in\sB(\mR^d\setminus\{0\}), t>0.
$$
The compensated Poisson random measure is given by $\tilde N((0,t]\times U)=N((0,t]\times U)-t\nu(U)$.
By L\'evy-It\^o's decomposition, one may write (cf. \cite{Sa})
\begin{align}
L_t=\int^t_0\!\!\!\int_{|x|\leq 1}x\tilde N(\dif s,\dif x)+\int^t_0\!\!\!\int_{|x|>1}xN(\dif s,\dif x).\label{Lev}
\end{align}
Let $X_t$ a purely discontinuous semimartingale with the form
\begin{align}
X_t=X_0+\int^t_0\xi_s\dif s+L_t,\label{Form}
\end{align}
where $X_0\in\sF_0$ and $(\xi_t)_{t\geq 0}$ is a measurable and $(\sF_t)$-adapted $\mR^d$-valued process.
Let $u$ be a  bounded smooth function on $\mR_+\times \mR^d$. By It\^o's formula (cf. \cite{Ap}), we have
\begin{align*}
u_t(X_t)&=u_0(X_0)+\int^t_0\Big([\p_s u_s+\cL_0 u_s](X_s)+\xi^i_s\p_i u_s(X_s)\Big)\dif s\\
&\quad+\int^t_0\!\!\!\int_{\mR^d\setminus\{0\}}[u_s(X_{s-}+y)-u_s(X_{s-})]\tilde N(\dif s,\dif y).
\end{align*}

In this section, we prove two estimates of Krylov's type for the above $X_t$.
Let us first prove the following simple Krylov's estimate, which will be used in Section 4
to prove the existence of weak solutions
for SDE (\ref{SDE}) with singular drift $b$.
\bt\label{Kry1}
Suppose that $\alpha\in(1,2)$,  $p>\frac{d}{\alpha-1}$ and $q>\frac{p\alpha}{p(\alpha-1)-d}$.
Then, for any $T_0>0$, there exist a constant $C=C(T_0,d,\alpha,p,q)>0$
such that for  any ($\sF_t$)-stopping time $\tau$, and $0\leq S<T\leq T_0$, and all $f\in L^q([S,T]; L^p(\mR^d))$,
\begin{align}
\mE\left(\int^{T\wedge\tau}_{S\wedge\tau} f_s(X_s)\dif s\Big|\sF_{S}\right)\leq
C\left(1+\mE\left(\int^{T\wedge\tau}_{S\wedge\tau} |\xi_s|\dif s\Big|\sF_S\right)\right)
\|f\|_{L^q([S,T];L^p(\mR^d))}.\label{Lp666}
\end{align}
\et
\begin{proof}
Let us first assume that $f\in C^\infty_0(\mR_+\times\mR^{d})$ and define
$$
u_t(x)=\int^t_0 \cT_{t-s}f_s(x)\dif s,
$$
where $\cT_t$ is defined by (\ref{Semi}).
By Lemma \ref{L11}, it is easy to see that $u_t(x)\in C^\infty(\mR_+\times\mR^d)$ and solves the following PIDE:
$$
\p_t u_t(x)=\cL_0u_t(x)+f_t(x).
$$
Choosing $\gamma\in(1+\frac{d}{p},\alpha-\frac{\alpha}{q})$,  by (\ref{Lp1}) and H\"older's inequality, we have
\begin{align}
\|u_t\|_{\gamma,p}&\leq \int^t_0  \|\cT_{t-s}f_s\|_{\gamma,p}\dif s
\leq C\int^t_0(t-s)^{-\gamma/\alpha}\|f_s\|_p\dif s\no\\
&\leq C\left(\int^t_0 (t-s)^{-q^*\gamma/\alpha}\dif s\right)^{\frac{1}{q^*}}
\|f\|_{L^q(\mR_+;L^p)}\leq C_t\|f\|_{L^q(\mR_+;L^p)},\label{Op01}
\end{align}
where $q^*=q/(q-1)$.

Fix $T_0>0$ and an ($\sF_t$)-stopping time $\tau$. Using It\^o's formula for $u_{T_0-t}(X_t)$
and by Doob's optimal theorem, we have
\begin{align*}
&\mE\Big(u_{T_0-T\wedge\tau}(X_{T\wedge\tau})|\sF_S\Big)-u_{T_0-S\wedge\tau}(X_{S\wedge\tau})\\
&\quad=\mE\left(\int^{T\wedge\tau}_{S\wedge\tau} \Big([\p_su_{T_0-s}+\cL_0u_{T_0-s}](X_s)
+\xi^i_s\p_i u_{T_0-s}(X_s)\Big)\dif s\Big|\sF_S\right)\\
&\quad\leq\mE\left(\int^{T\wedge\tau}_{S\wedge\tau} \Big(-f_s(X_s)+|\xi_s|\cdot|\nabla u_{T_0-s}|(X_s)\Big)\dif s\Big|\sF_S\right),
\end{align*}
which yields by (\ref{Op01}) and (\ref{em}) that,
\begin{align*}
\mE\left(\int^{T\wedge\tau}_{S\wedge\tau}  f_s(X_s)\dif s\Big|\sF_S\right)&\leq 2\sup_{s,x}|u_s(x)|
+\sup_{s,x}|\nabla u_s|(x)\cdot\mE\left(\int^{T\wedge\tau}_{S\wedge\tau} |\xi_s|\dif s\Big|\sF_S\right)\\
&\leq C\|f\|_{L^q(\mR_+;L^p)} \left(1+\mE\left(\int^{T\wedge\tau}_{S\wedge\tau} |\xi_s|\dif s\Big|\sF_S\right)\right),
\end{align*}
where we have used  $p(\gamma-1)>d$.
By a standard density argument, we obtain (\ref{Lp666}) for general $f\in L^q([S,T]; L^p(\mR^d))$.
\end{proof}

In one dimensional case, as in \cite{Ku1}, we even have:
\bt
Let $X_t$ take the following form:
\begin{align*}
X_t=X_0+\int^t_0\xi_s\dif s+\int^t_0h_s\dif L_s,
\end{align*}
where $h_s$ is a bounded predictable process. Suppose that $\alpha\in(1,2)$ and $p>\frac{1}{\alpha-1}$.
Then, for any $T_0>0$, there exist a constant $C=C(T_0,\alpha,p,q)>0$
such that for any ($\sF_t$)-stopping time $\tau$, and $0\leq S<T\leq T_0$, and all $f\in L^p(\mR^d)$,
\begin{align}
\mE\left(\int^{T\wedge\tau}_{S\wedge\tau}|h_s|^\alpha f(X_s)\dif s\Big|\sF_S\right)\leq
C\left(1+\mE\left(\int^{T\wedge\tau}_{S\wedge\tau}(|\xi_s|+|h_s|^\alpha)\dif s\Big|\sF_S\right)\right)\|f\|_p.\label{Lp6666}
\end{align}
\et
\begin{proof}
Fix $T_0>0$. For $f\in C^\infty_0(\mR^{d})$, let us define
$$
u_{T_0}(x):=\int^{T_0}_0\cT_{T_0-s}f(x)\dif s.
$$
It is easy to see that
\begin{align}
\cL_0 u_{T_0}(x)=\cT_{T_0}f(x)-f(x).\label{PL1}
\end{align}
Using It\^o's formula for $u_{T_0}(X_t)$ (cf. \cite[Proposition 2.1]{Ba1}), one finds that
$$
\mE \left(u_{T_0}(X_{t\wedge\tau})|\sF_S\right)=u_{T_0}(X_{S\wedge\tau})+\mE\left(\int^{T\wedge\tau}_{S\wedge\tau}
(|h_s|^\alpha\cL_0 u_{T_0}(X_s)+\xi_s u'_{T_0}(X_s))\dif s\Big|\sF_S\right),
$$
which together with (\ref{PL1}) yields that
\begin{align*}
\mE\left(\int^{T\wedge\tau}_{S\wedge\tau}|h_s|^\alpha f(X_s)\dif s\Big|\sF_S\right)&\leq
2\|u_{T_0}\|_\infty+\|u'_{T_0}\|_\infty\mE\left(\int^{T\wedge\tau}_{S\wedge\tau}\xi_s \dif s\Big|\sF_S\right)\\
&\quad+\|\cT_{T_0}f\|_\infty\mE\left(\int^{T\wedge\tau}_{S\wedge\tau}|h_s|^\alpha\dif s\Big|\sF_S\right)\\
&\leq C\left(1+\mE\left(\int^{T\wedge\tau}_{S\wedge\tau}(|\xi_s|+|h_s|^\alpha)\dif s\Big|\sF_S\right)\right)\|f\|_p,
\end{align*}
where we have used $p(\alpha-1)>1$, (\ref{Lp1}) and (\ref{em}).
By a standard density argument, we obtain (\ref{Lp6666}) for general $f\in L^p(\mR^d)$.
\end{proof}
In the above two theorems, the requirement of $p>\frac{d}{\gamma-1}$ is too strong to prove Theorem \ref{Main}.
It is clear that this is caused by directly controlling the $\infty$-norm of $\nabla u_s(x)$
by Sobolev embedding theorem. In what follows, we shall relax it to $p>\frac{d}{\gamma}$. The price to pay is that
we need to assume that $\xi_t$ is a bounded ($\sF_t$)-adapted process. Nevertheless,
Theorem \ref{Kry1} can be used to prove the existence of weak solutions for
SDE (\ref{SDE}) with globally integrable drift.

We now start by solving the following semi-linear PIDE:
\begin{align}
\p_t u=\cL_0 u+\kappa|\nabla u|+f,\ \ u_0\equiv 0,\ \ t\geq 0\label{PP1}
\end{align}
where $\kappa>0$, $\cL_0$ is the  generator of L\'evy process $(L_t)_{t\geq 0}$ given by (\ref{Op4}),
and $f$ is a locally integrable function on $\mR_+\times\mR^d$.

We first give the following definition of generalized solutions to PIDE (\ref{PP1}).
\bd
For $p\geq 1$, a function $u\in C([0,\infty);\mH^1_p)$ is called a generalized solution of (\ref{PP1}),
if for all function $\varphi\in C^\infty_0([0,\infty)\times\mR^d)$, it holds that
$$
-\int^\infty_0\!\!\!\int_{\mR^d}u\p_t\varphi=\int^\infty_0\!\!\!\int_{\mR^d}u\cL_0^*\varphi
+\int^\infty_0\!\!\!\int_{\mR^d}(\kappa|\nabla u|+f)\varphi,
$$
where $\cL^*_0$ is the adjoint operator of $\cL_0$ given by
$$
\cL^*_0\varphi(x)=\int_{\mR^d}(\varphi(x-z)-\varphi(x)+1_{|z|\leq 1}z^i\p_i\varphi(x))\nu(\dif z).
$$
\ed
\br\label{Re}
If we extend $u$ and $f$ to $\mR$ by setting $u_t=f_t\equiv 0$ for $t\leq 0$, then
for all $\varphi\in C^\infty_0(\mR^{d+1})$,
$$
-\int_{\mR^{d+1}}u\p_t\varphi=\int_{\mR^{d+1}}u\cL_0^*\varphi
+\int_{\mR^{d+1}}(\kappa|\nabla u|+f)\varphi.
$$
Since the L\'evy measure $\nu$ is symmetric, $\cL^*_0$ is in fact the same as $\cL_0$.
\er
The following proposition is now standard. We omit the proof.
\bp
For $p\geq 1$, let $u\in C([0,\infty);\mH^1_p)$ and $f\in L^1_{loc}([0,\infty)\times\mR^d)$. The following three statements are equivalent:
\begin{enumerate}[(i)]
\item $u$ is a generalized solution of (\ref{PP1});
\item For any $\phi\in C^\infty_0(\mR^d)$, it holds that for all $t\geq 0$,
$$
\int_{\mR^d}u_t\phi=\int^t_0\!\!\!\int_{\mR^d}u_s\cL_0^*\phi
+\int^t_0\!\!\!\int_{\mR^d}(\kappa|\nabla u_s|+f)\phi;
$$
\item $u$ satisfies the following integral equation:
$$
u_t(x)=\int^t_0\cT_{t-s}(\kappa|\nabla u_s|+f_s)(x)\dif s,\ \ \forall t\geq 0.
$$
\end{enumerate}
\ep

We have the following existence-uniqueness result about the generalized solution of equation (\ref{PP1}).
\bt
For $p>1$, $\alpha\in(1,2)$, $\gamma\in[1,\alpha)$ and $q>\frac{\alpha}{\alpha-\gamma}$, assume
that $f\in L^q_{loc}(\mR_+; L^p(\mR^d))$. Then,
there exists a unique generalized solution $u\in C([0,\infty);\mH^\gamma_p)$ to PIDE (\ref{PP1}). Moreover,
\begin{align}
\|u_t\|_{\gamma,p}\leq C_t\|f\|_{L^q([0,t]; L^p)},\ \ \forall t\geq 0,\label{PP4}
\end{align}
where $C_t\geq 0$ is a continuous increasing function of $t$ with $C_t=O(t^{1-\frac{\gamma}{\alpha}-\frac{1}{q}})$ as $t\to 0$.
\et
\begin{proof}
Let $u^{(0)}\equiv 0$. For $n\in\mN$, define $u^{(n)}$ recursively by
\begin{align}
u^{(n)}_t(x)=\int^t_0\cT_{t-s}(\kappa|\nabla u^{(n-1)}_s|+f_s)(x)\dif s,\ \ \forall t\geq 0.\label{PP3}
\end{align}
By (i) of Lemma \ref{L11} and H\"older's inequality, we have for $q>\frac{\alpha}{\alpha-\gamma}$,
\begin{align*}
\|u^{(n)}_t\|_{\gamma,p}&\leq C\int^t_0(t-s)^{-\frac{\gamma}{\alpha}}\Big(\kappa\|\nabla u^{(n-1)}_s\|_p+\|f_s\|_p\Big)\dif s\\
&\leq C\left(\int^t_0(t-s)^{-\frac{q\gamma}{(q-1)\alpha}}\dif s\right)^{\frac{q-1}{q}}
\left(\int^t_0\Big(\kappa^q\|\nabla u^{(n-1)}_s\|^q_p+\|f_s\|^q_p\Big)\dif s\right)^{\frac{1}{q}},
\end{align*}
which yields that
\begin{align*}
\|u^{(n)}_t\|_{\gamma,p}^q&\leq Ct^{\frac{q(\alpha-\gamma)-\alpha}{\alpha}}
\left(\int^t_0\|u^{(n-1)}_s\|^q_{1,p}\dif s+\int^t_0\|f_s\|^q_p\dif s\right)\\
&\leq Ct^{\frac{q(\alpha-\gamma)-\alpha}{\alpha}}\left(\int^t_0\|u^{(n-1)}_s\|^q_{\gamma,p}\dif s
+\int^t_0\|f_s\|^q_p\dif s\right).
\end{align*}
By Gronwall's inequality, we obtain that for all $t\geq 0$,
\begin{align}
\sup_{n\in\mN}\|u^{(n)}_t\|_{\gamma,p}^q\leq C_t\int^t_0\|f_s\|^q_p\dif s.\label{PP2}
\end{align}

Next, fixing $T>0$, we want to prove the H\"older continuity of mapping $[0,T]\ni t\mapsto u_t^{(n)}\in\mH^\gamma_p$.
For $T\geq t>t'\geq 0$, we have
\begin{align*}
u^{(n)}_t-u^{(n)}_{t'}&=\int^{t'}_0(\cT_{t-s}-\cT_{t'-s})(\kappa|\nabla u^{(n-1)}_s|+f_s)\dif s\\
&\quad+\int^t_{t'}\cT_{t-s}(\kappa|\nabla u^{(n-1)}_s|+f_s)\dif s=:I_1(t,t')+I_2(t,t').
\end{align*}
For $I_1(t,t')$, using the semigroup property of $\cT_t$, we further have
$$
I_1(t,t')=\int^{t'}_0\cT_{(t'-s)/2}(\cT_{t-t'}-I)\cT_{(t'-s)/2}(\kappa|\nabla u^{(n-1)}_s|+f_s)\dif s.
$$
Hence, by Lemma \ref{L11} and (\ref{PP2}), for $\delta\in(0,\alpha-\gamma-\frac{\alpha}{q})$, we have
\begin{align}
\|I_1(t,t')\|_{\gamma,p}&\leq C\int^{t'}_0(t'-s)^{-\frac{\gamma}{\alpha}}
\|(\cT_{t-t'}-I)\cT_{(t'-s)/2}(\kappa|\nabla u^{(n-1)}_s|+f_s)\|_p\dif s\no\\
&\leq C\int^{t'}_0(t'-s)^{-\frac{\gamma}{\alpha}}(t-t')^{\frac{\delta}{\alpha}}
\|\cT_{(t'-s)/2}(\kappa|\nabla u^{(n-1)}_s|+f_s)\|_{\delta,p}\dif s\no\\
&\leq C(t-t')^{\frac{\delta}{\alpha}}\int^{t'}_0(t'-s)^{-\frac{\gamma+\delta}{\alpha}}
(\|\nabla u^{(n-1)}_s\|_p+\|f_s\|_p)\dif s\no\\
&\leq C(t-t')^{\frac{\delta}{\alpha}}\|f\|_{L^q([0,T];L^p(\mR^d))}.\label{Lp3}
\end{align}
For $I_2(t,t')$, using (\ref{PP2}), we also have
\begin{align}
\|I_2(t,t')\|_{\gamma,p}\leq C_T(t-t')^{1-\frac{\gamma}{\alpha}-\frac{1}{q}}\|f\|_{L^q([0,T]; L^p)}.\label{Lp4}
\end{align}
Combining (\ref{Lp3}) and (\ref{Lp4}), we obtain the desired H\"older continuity.

Now, as above, we can make the following estimation:
\begin{align*}
\|u^{(n)}_t-u^{(m)}_t\|_{\gamma,p}&\leq C\int^t_0(t-s)^{-\frac{\gamma}{\alpha}}
\Big(\||\nabla u^{(n-1)}_s|-|\nabla u^{(m-1)}_s|\|_p\Big)\dif s\\
&\leq Ct^{1-\frac{\gamma}{\alpha}-\frac{1}{q}}\left(\int^t_0\|\nabla (u^{(n-1)}_s-u^{(m-1)}_s)\|^q_p\dif s\right)^{\frac{1}{q}},
\end{align*}
which then gives that
$$
\|u^{(n)}_t-u^{(m)}_t\|^q_{\gamma,p}\leq Ct^{\frac{q(\alpha-\gamma)-\alpha}{\alpha}}
\int^t_0\|u^{(n-1)}_s-u^{(m-1)}_s\|^q_{\gamma,p}\dif s,
$$
where $C$ is independent of $n,m$ and $t$.
Using (\ref{PP2}) and Fatou's lemma, we find that
\begin{align*}
\varlimsup_{n,m\to\infty}\sup_{s\in[0,t]}\|u^{(n)}_s-u^{(m)}_s\|^q_{\gamma,p}&\leq Ct^{\frac{q(\alpha-\gamma)-\alpha}{\alpha}}\int^t_0
\varlimsup_{n,m\to\infty}\|u^{(n-1)}_s-u^{(m-1)}_s\|^q_{\gamma,p}\dif s\\
&\leq Ct^{\frac{q(\alpha-\gamma)-\alpha}{\alpha}}\int^t_0\varlimsup_{n,m\to\infty}\sup_{r\in[0,s]}\|u^{(n-1)}_r-u^{(m-1)}_r\|^q_{\gamma,p}\dif s,
\end{align*}
and so, for any $t>0$,
$$
\varlimsup_{n,m\to\infty}\sup_{s\in[0,t]}\|u^{(n)}_s-u^{(m)}_s\|^q_{\gamma,p}=0.
$$
Thus, there exists a $u\in C([0,\infty);\mH^\gamma_p)$ such that for any $t>0$,
$$
\lim_{n\to\infty}\sup_{s\in[0,t]}\|u^{(n)}_s-u_s\|_{\gamma,p}=0.
$$
Taking limits for both sides of (\ref{PP3}), we obtain the existence of a generalized solution,
and (\ref{PP4}) is direct from (\ref{PP2}).

As for the uniqueness, it follows from a similar calculation.
The proof is complete.
\end{proof}

Let us now prove our second Krylov's estimate.
\bt\label{Th1}
Suppose that $\alpha\in(1,2)$, $p>\frac{d}{\alpha}\vee 1$ and $q>\frac{p\alpha}{p\alpha-d}$.
Let $(\xi_t)_{t\geq 0}$ be a measurable and ($\sF_t$)-adapted process bounded by $\kappa$, and
let $X_t$ have the form (\ref{Form}). Then for any $T_0>0$, there exist a constant $C=C(T_0,\kappa, d,\alpha,p,q)>0$
such that for  any ($\sF_t$)-stopping time $\tau$, and $0\leq S<T\leq T_0$, and all $f\in L^q([S,T]; L^p(\mR^d))$,
\begin{align}
\mE\left(\int^{T\wedge\tau}_{S\wedge\tau} f_s(X_s)\dif s\Big|\sF_S\right)\leq C\|f\|_{L^q([S,T];L^p(\mR^d))}.\label{Lp66}
\end{align}
\et
\begin{proof}
Let us first assume that $f\in C^\infty_0(\mR_+\times\mR^{d})$. Choose $\gamma\in(\frac{d}{p},\alpha-\frac{\alpha}{q})$
and let $u\in C([0,\infty);\mH^\gamma_p)$ be the unique solution of  PIDE (\ref{PP1}).
Fix $T_0>0$, and let $v_t(x)=u_{T_0-t}(x)$.
It is easy to see that $v_t$ is a generalized solution of the following PIDE:
\begin{align}
\p_tv+\cL_0 v+\kappa|\nabla v|+f=0,\ \ v_{T_0}\equiv 0.\label{Lp5}
\end{align}
Let $\rho$ be a smooth nonnegative function in $\mR^{d+1}$ with support in
$\{(s,x)\in\mR^{d+1}: |s|+|x|\leq 1\}$ and $\int_{\mR^{d+1}}\rho=1$. For $\eps>0$, set
$$
\rho_\eps(s,x)=\eps^{-(d+1)}\rho(\eps^{-1}s,\eps^{-1}x)
$$
and
$$
v^{(\eps)}=v*\rho_\eps,\ \ f^{(\eps)}=f*\rho_\eps.
$$
Taking convolutions for both sides of (\ref{Lp5}), we obtain that
$$
\p_t v^{(\eps)}+\cL_0v^{(\eps)}+\kappa |\nabla v^{(\eps)}|+f^{(\eps)}\leq (\p_tv+\cL_0 v+\kappa|\nabla v|+f)*\rho_\eps=0.
$$
Here we have used Remark \ref{Re}.

Using It\^o's formula for $v^{(\eps)}_t(X_t)$, we get
\begin{align*}
\mE (v^{(\eps)}_{T\wedge\tau}(X_{T\wedge\tau})|\sF_S)-v^{(\eps}_{S\wedge\tau}(X_{S\wedge\tau})
&=\mE\left(\int^{T\wedge\tau}_{S\wedge\tau}\Big( [\p_sv^{(\eps)}_s+\cL_0v^{(\eps)}_s](X_S)+\xi^i_s\p_i v^{(\eps)}_s(X_s)\Big)\dif s\Big|\sF_S\right)\\
&\leq\mE\left(\int^{T\wedge\tau}_{S\wedge\tau}[\p_sv^{(\eps)}_s+\cL_0v^{(\eps)}_s+\kappa |\nabla v^{(\eps)}_s|](X_s)\dif s\Big|\sF_S\right)\\
&\leq-\mE\left(\int^{T\wedge\tau}_{S\wedge\tau} f^{(\eps)}_s(X_s)\dif s\Big|\sF_S\right),
\end{align*}
which yields by (\ref{PP4}) and (\ref{em}) that,
\begin{align*}
\mE\left(\int^{T\wedge\tau}_{S\wedge\tau} f^{(\eps)}_s(X_s)\dif s\Big|\sF_S\right)&\leq
2\sup_{(s,x)\in[0,T_0]\times\mR^d}|v^{(\eps)}_s(x)|\leq
2\sup_{(t,x)\in[0,T_0]\times\mR^d}|v_t(x)|\leq\\
&\leq2\sup_{(t,x)\in[0,T_0]\times\mR^d}|u_t(x)|\leq C\int^{T_0}_0\|f_s\|^q_{L^p}\dif s.
\end{align*}
Taking limits $\eps\to 0$, by the dominated convergence theorem, we have
$$
\mE\left(\int^{T\wedge\tau}_{S\wedge\tau} f_s(X_s)\dif s\Big|\sF_S\right)\leq C\int^{T_0}_0\|f_s\|^q_{L^p}\dif s.
$$
By a standard density argument, we obtain (\ref{Lp66}) for general $f\in L^q([S,T]; L^p(\mR^d))$.
\end{proof}

\section{Weak solutions for SDE (\ref{SDE}) with globally integrable drift}

In this section, we use Theorem \ref{Kry1} to prove the following existence of weak solutions for  SDE (\ref{SDE}).
\bt\label{Weak}
Suppose that $\alpha\in(1,2)$, $\gamma\in(1,\alpha)$, $p>\frac{d}{\gamma-1}$ and $q>\frac{\alpha}{\alpha-\gamma}$.
Then for any  $b\in L^\infty_{loc}(\mR_+; L^\infty(\mR^d))+L^q_{loc}(\mR_+;L^p(\mR^d))$ and $x_0\in\mR^d$, there
exists a weak solution to SDE (\ref{SDE}). More precisely, there exists a probability space $(\tilde\Omega,\tilde\sF,\tilde P)$
and two c\`adl\`ag stochastic processes $\tilde X_t$ and $\tilde L_t$ defined on it such that
$\tilde L_t$ is a symmetric $\alpha$-stable process with respect to the completed filtration
$\tilde\sF_t:=\sigma^{\tilde P}\{\tilde X_s,\tilde L_s, s\leq t\}$ and
$$
\tilde X_t=x_0+\int^t_0b(s,\tilde X_s)\dif s+\tilde L_t\ \ \forall t\geq 0.
$$
\et
\begin{proof}
Our proof is adapted from the proof of \cite[p.87, Theorem 1]{Kr}.
Let $b=b_1+b_2$ with $b_1\in L^\infty_{loc}(\mR_+; L^\infty(\mR^d))$ and
$b_2\in L^q_{loc}(\mR_+;L^p(\mR^d))$.
Let $b^{(n)}_i(t,x)=(b_i(t,\cdot)*\rho_n)(x)$ be the mollifying approximation of $b_i$, $i=1,2$.
It is easy to see that for some $\ell^{(n)}_t\in L^1_{loc}(\mR_+)$,
$$
|b^{(n)}(t,x)-b^{(n)}(t,x)|\leq \ell_t|x-y|,\ \ \forall x,y\in\mR^d.
$$
Let $X^{(n)}_t$ solve the following SDE:
$$
X^{(n)}_t=x_0+\int^t_0b^{(n)}(s,X^{(n)}_s)\dif s+L_t.
$$

{\bf (Claim 1:)} For some $\delta>1$, we have
\begin{align}
\sup_{n\in\mN}\mE\int^T_0|b^{(n)}(s,X^{(n)}_s)|^\delta\dif s<+\infty,\ \ \forall T>0.\label{Es4}
\end{align}
In fact, choosing $\delta>1$ and $p'\in(\frac{d}{\gamma-1},p), q'\in(\frac{\alpha}{\alpha-\gamma},q)$ such that
$p'\delta=p$ and $q'\delta=q$, by (\ref{Lp666}) and Young's inequality, we have
\begin{align*}
\mE\int^T_0|b^{(n)}_2(s,X^{(n)}_s)|^\delta\dif s&\leq C_T
\left(1+\mE\int^T_0|b^{(n)}(s,X^{(n)}_s)|\dif s\right)\||b^{(n)}_2|^\delta\|_{L^{q'}([0,T];L^{p'}(\mR^d))}\no\\
&\leq C_T\left(1+\|b_1\|_{L^\infty([0,T]\times\mR^d)}+\mE\int^T_0|b^{(n)}_2(s,X^{(n)}_s)|\dif s\right)
\|b_2\|^\delta_{L^q([0,T];L^p(\mR^d))}\\
&\leq \frac{1}{2}\mE\int^T_0|b^{(n)}_2(s,X^{(n)}_s)|^\delta\dif s
+C_T\|b_2\|^{\frac{\delta^2}{\delta-1}}_{L^q([0,T];L^p(\mR^d))}\\
&\quad+C_T\left(1+\|b_1\|_{L^\infty([0,T]\times\mR^d)}\right)\|b_2\|^\delta_{L^q([0,T];L^p(\mR^d))},
\end{align*}
which then implies (\ref{Es4}).

Let $\mD$ be the space of all c\`adl\`ag functions from $\mR_+$ to $\mR^d$,
which is endowed with the Skorohod topology so that $\mD$ is a Polish space.
Set
$$
H^{(n)}_t:=\int^t_0b^{(n)}(s,X^{(n)}_s)\dif s.
$$
Using Claim 1, it is easy to check that the following
Aldous'  tightness criterions \cite{Ad} hold:
$$
\lim_{N\to\infty}\varlimsup_{n\to\infty} P\left(\sup_{t\in[0,T]}|H^{(n)}_t|\geq N\right)=0, \ \forall T>0,
$$
and
$$
\lim_{\eps\to 0}\varlimsup_{n\to\infty} \sup_{\tau\in\cS_T}
P\left(|H^{(n)}_{\tau}-|H^{(n)}_{\tau+\eps}|\geq a\right)=0,\ \ \forall T,a>0,
$$
where $\cS_T$ denotes all the bounded stopping times with bound $T$.
Thus, the law of $t\mapsto H^{(n)}_t$ in $\mD$ is tight, and so does $(H^{(n)}_\cdot, L_\cdot)$.
By Prohorov's theorem, there exists a subsequence still denoted by $n$ such that the law of $(H^{(n)}_\cdot, L_\cdot)$
in $\mD\times\mD$ weakly converges, which then implies that the law of $(X^{(n)}_\cdot, L_\cdot)$
weakly converges. By Skorohod's representation theorem, there is a probability space $(\tilde\Omega,\tilde\sF,\tilde P)$
and the $\mD\times\mD$-valued random variables $(\tilde X^{(n)}_\cdot, \tilde L^{(n)}_\cdot)$ and
$(\tilde X_\cdot, \tilde L_\cdot)$ such that

(i) $(\tilde X^{(n)}_\cdot, \tilde L^{(n)}_\cdot)$ has the same law as $(X^{(n)}_\cdot, L_\cdot)$ in $\mD\times\mD$;

(ii) $(\tilde X^{(n)}_\cdot, \tilde L^{(n)}_\cdot)$ converges to
$(\tilde X_\cdot, \tilde L_\cdot)$, $\tilde P$-almost surely.

In particular, $\tilde L$ is still a symmetric $\alpha$-stable process and
$$
\tilde X^{(n)}_t=x_0+\int^t_0b^{(n)}(s,\tilde X^{(n)}_s)\dif s+\tilde L^{(n)}_t.
$$

{\bf (Claim 2:)} For any nonnegative measurable function $f$ and $T>0$, we have
$$
\tilde\mE\int^T_0 f_s(\tilde X_s)\dif s\leq C_T \|f\|_{L^q([0,T];L^p(\mR^d))},
$$
where $\tilde\mE$ denotes the expectation with respect to the probability measure $\tilde P$.

Let $f\in C_0([0,T]\times\mR^d)$. By the dominated convergence theorem, we have
\begin{align*}
\tilde\mE\int^T_0 f_s(\tilde X_s)\dif s
&=\lim_{n\to\infty}\tilde\mE\int^T_0 f_s(\tilde X^{(n)}_s)\dif s\\
&=\lim_{n\to\infty}\mE\int^T_0 f_s(X^{(n)}_s)\dif s\\
&\leq C \|f\|_{L^q([0,T];L^p(\mR^d))},
\end{align*}
where in the last step we have used (\ref{Lp666}) and (\ref{Es4}).
For general $f$, it follows by the monotone class theorem.

The proof will be finished if one can show the following claim:

{\bf (Claim 3:)} For any $T>0$, we have
\begin{align}
\lim_{n\to\infty}\tilde \mE\left(\int^T_0 |b^{(n)}_i(s,\tilde X^{(n)}_s)
-b_i(s,\tilde X_s)|\dif s\right)=0,\ \ i=1,2.\label{Lo1}
\end{align}
Let $\chi_R(x)$ be a smooth nonnegative function on $\mR^d$ with $\chi_R(x)=1$ for $|x|\leq R$ and
$\chi_R(x)=0$ for $|x|>R+1$. Then for any $n,m\in\mN$,
\begin{align}
\tilde \mE\left(\int^T_0 |b^{(n)}_1(s,\tilde X^{(n)}_s)
-b_1(s,\tilde X_s)|\dif s\right)&\leq\tilde\mE\left(\int^T_0|b^{(n)}_1(s,\tilde X^{(n)}_s)
-b^{(m)}_1(s,\tilde X^{(n)}_s)|\dif s\right)\no\\
&\quad+\tilde \mE\left(\int^T_0|b^{(m)}_1(s,\tilde X^{(n)}_s)
-b^{(m)}_1(s,\tilde X_s)|\dif s\right)\no\\
&\quad+\tilde \mE\left(\int^T_0|b^{(m)}_1(s,\tilde X_s)
-b_1(s,\tilde X_s)|\dif s\right)\no\\
&=:I^{(n,m)}_1+I^{(n,m)}_2+I^{(n,m)}_3.\label{Lp11}
\end{align}
For fixed $m$, by the above (ii) and the dominated convergence theorem, we have
$$
\lim_{n\to\infty}I^{(n,m)}_2=0.
$$
For $I^{(n,m)}_1$, by Claim 1, we have
\begin{align*}
I^{(n,m)}_1&\leq \|b_1\|_{L^\infty([0,T];L^\infty(\mR^d)}\tilde \mE\left(\int^T_0|1-\chi_R(\tilde X^{(n)}_s)|\dif s\right)+
\tilde \mE\left(\int^T_0[\chi_R(|b^{(n)}_1-b^{(m)}_1)](s,\tilde X^{(n)}_s)|\dif s\right)\\
&\leq \frac{C}{R}\int^T_0 \mE|X^{(n)}_s|\dif s+
\mE\left(\int^T_0[\chi_R(|b^{(n)}_1-b^{(m)}_1)](s,X^{(n)}_s)|\dif s\right)\\
&\leq \frac{C}{R}+C\|\chi_R(|b^{(n)}_1-b^{(m)}_1)\|_{L^q([0,T];L^p)}.
\end{align*}
Similarly, by Claim 2, we have
$$
I^{(n,m)}_3\leq\frac{C}{R}+C\|\chi_R(|b^{(m)}_1-b_1)\|_{L^q([0,T];L^p)}.
$$
Taking limits for both sides of (\ref{Lp11}) in order:
$n\to\infty$, $m\to\infty$ and $R\to\infty$, we obtain (\ref{Lo1}) for $i=1$.
It is similar to prove (\ref{Lo1}) for $i=2$. The whole proof is complete.
\end{proof}

\br
When $b$ is time-independent and the L\'evy measure $\nu(\dif\xi)=\frac{C_\alpha}{|\xi|^{d+\alpha}}\dif\xi$,
Theorem \ref{Weak} has been proven by Chen, Kim and Song \cite[Theorem 2.5]{Ch-Ki-So} by different argument.
\er
\section{Proof of Theorem \ref{Main}}

We now consider the following linear PIDE for $\lambda>0$:
\begin{align}
\p_t u=(\cL_0-\lambda) u+b^i\p_iu+f,\ \ u_0\equiv0.\label{PP6}
\end{align}
As in the previous section, one may define the notion of generalized solutions and has:
\bt
Let $\alpha\in(1,2)$ and $\gamma\in(1,\alpha)$. Assume that for some $p>\frac{d}{\gamma}$
and $0\leq\beta\in(1-\gamma+\frac{d}{p},1)$,
$$
b\in L^\infty_{loc}(\mR_+,L^\infty(\mR^d)\cap\mW^\beta_p),\ \ f\in L^\infty_{loc}(\mR_+; \mW^\beta_p).
$$
Then, there exists a unique generalized
solution $u=u^\lambda\in C(\mR_+;\mH^{\gamma+\beta}_p)$ to PIDE (\ref{PP6}).
Moreover, for some $\delta>0$ and any $\lambda\geq 1$,
\begin{align}
\|u^\lambda_t\|_{\gamma+\beta,p}\leq C_t\lambda^{-\delta}
\|f\|_{L^\infty([0,t]; \mH^\beta_p)},\ \ \forall t\geq 0,\label{PP5}
\end{align}
where $C_t>0$ is an increasing function of $t$ with $\lim_{t\downarrow 0}C_t=0$.
\et
\begin{proof}
As in the proof of Theorem \ref{Th1}, we only need to prove the a priori estimate (\ref{PP5}).
Let $u$ satisfy the following integral equation:
$$
u_t(x)=\int^t_0e^{-\lambda(t-s)}\cT_{t-s}(b^i_s\p_iu_s+f_s)(x)\dif s,\ \ \forall t\geq 0.
$$
Let $\eps\in(0,\alpha-\gamma)$ and $q>\frac{\alpha}{\alpha-\gamma-\eps}$.
By Lemma \ref{L11} and H\"older's inequality, we have
\begin{align*}
\|u_t\|_{\gamma+\beta,p}&\leq C\int^t_0e^{-\lambda(t-s)}(t-s)^{-\frac{\gamma+\eps}{\alpha}}
\Big(\|b^i_s\p_iu_s\|_{\beta-\eps,p}+\|f_s\|_{\beta-\eps,p}\Big)\dif s\\
&\stackrel{(\ref{Fr})}{\leq} C\left(\int^t_0e^{-\lambda q(t-s)}(t-s)^{-\frac{q(\gamma+\eps)}{(q-1)\alpha}}\dif s\right)^{\frac{q-1}{q}}
\left(\int^t_0\Big(\|b^i_s\p_iu_s\|^\sim_{\beta,p}+\|f_s\|^\sim_{\beta,p}\Big)^q\dif s\right)^{\frac{1}{q}}\\
&\leq C\lambda^{\frac{\gamma+\eps}{\alpha}-1+\frac{1}{q}}\left(\int^\infty_0e^{-s}
s^{-\frac{q(\gamma+\eps)}{(q-1)\alpha}}\dif s\right)^{\frac{q-1}{q}}
\left(\int^t_0\Big(\|b^i_s\p_iu_s\|_{\beta,p}^\sim+\|f_s\|_{\beta,p}^\sim\Big)^q\dif s\right)^{\frac{1}{q}}.
\end{align*}
In view of $(\gamma+\beta-1)p>d$ and $\gamma>1$, we have
\begin{align*}
\|b^i_s\p_iu_s\|_{\beta,p}^\sim&\stackrel{(\ref{LL1})}{\leq}
\|b_s\|_\infty\|\nabla u_s\|_p+\left(\int_{\mR^d}\!\!\int_{\mR^d}
\frac{|(b^i_s\p_iu_s)(x)-(b^i_s\p_iu_s)(y)|^p}{|x-y|^{d+\beta p}}\dif x\dif y\right)^{1/p}\\
&\leq \|b_s\|_\infty\|\nabla u_s\|_p+\|b_s\|_\infty\|\nabla u_s\|^\sim_{\beta,p}
+\|b_s\|^\sim_{\beta,p}\|\nabla u_s\|_\infty\\
&\stackrel{(\ref{em})}{\leq}  \|b_s\|_\infty\|u_s\|_{1,p}+\|b_s\|_\infty\|u_s\|^\sim_{1+\beta,p}
+C\|b_s\|^\sim_{\beta,p}\|u_s\|_{\gamma+\beta,p}\\
&\stackrel{(\ref{Fr})}{\leq}  C(\|b_s\|_\infty+\|b_s\|^\sim_{\beta,p})\|u_s\|_{\gamma+\beta,p}.
\end{align*}
Hence,
\begin{align*}
\|u_t\|_{\gamma+\beta,p}^q
&\leq C\lambda^{\frac{q(\gamma+\eps)}{\alpha}-q+1}\left(
\|b\|^q_{L^\infty([0,t];L^\infty\cap\mW^\beta_p)}\int^t_0\|u_s\|^q_{\gamma+\beta,p}\dif s
+t\|f\|^q_{L^\infty([0,t];\mW^\beta_p)}\right).
\end{align*}
By Gronwall's inequality, we obtain (\ref{PP5}) with $\delta=q-1-\frac{q(\gamma+\eps)}{\alpha}>0$.
\end{proof}

Below, we assume that $b\in L^\infty_{loc}(\mR_+;L^\infty(\mR^d)\cap\mW^\beta_p)$ with
\begin{align}
\beta\in(1-\frac{\alpha}{2},1),\ \ \ p>\frac{2d}{\alpha},\label{Con}
\end{align}
and fix
$$
\gamma\in\Big((1+\frac{\alpha}{2}-\beta)\vee 1,\alpha\Big).
$$
Let $u^{\ell}$ solve the following PIDE:
\begin{align*}
\p_t u^{\ell}=(\cL_0-\lambda) u^{\ell}+b^i\p_iu^{\ell}+b^\ell,\ \ u^{\ell}_0(x)=0,\ \ \ell=1,\cdots, d.
\end{align*}
Fix $T>0$ and set
$$
\v_t(x):=(u^1_{T-t}(x),\cdots,u^d_{T-t}(x)).
$$
Then $\v_t(x)$ solves the following PIDE:
\begin{align}
\p_t\v+(\cL_0-\lambda)\v+b^i\p_i\v+b=0,\ \ \v_T(x)=0.\label{PP7}
\end{align}
Since $(\gamma+\beta-1)p>d$, by (\ref{em}) and (\ref{PP5}),
one can choose $\lambda$ sufficiently large such that
\begin{align}
\sup_{t\in[0,T]}\sup_{x\in\mR^d}|\nabla\v_t(x)|\leq \sup_{t\in[0,T]}C\|\v_t\|_{\gamma+\beta,p}
\leq C_T\lambda^{-\delta}\|f\|_{L^\infty([0,T]; \mH^\beta_p)}\leq\frac{1}{2}.\label{Lp8}
\end{align}
Let us define
$$
\Phi_t(x)=x+\v_t(x).
$$
Since for each $t\in[0,T]$,
$$
\frac{1}{2}|x-y|\leq|\Phi_t(x)-\Phi_t(y)|\leq\frac{3}{2}|x-y|,
$$
$x\mapsto\Phi_t(x)$ is a diffeomorphism and
\begin{align}
|\nabla\Phi_t(x)|\leq\frac{3}{2},\ \ |\nabla\Phi^{-1}_t(x)|\leq 2.\label{Op1}
\end{align}
\bl\label{Le3}
Let $\Phi_t(x)$ be defined as above. Fix an ($\sF_t$)-stopping time $\tau$ and let $X_t\in\sS_b^\tau(x)$ be
a local solution of SDE (\ref{SDE}). Then $Y_t=\Phi_t(X_t)$ solves the following SDE on $[0,T\wedge\tau)$:
\begin{align}
Y_t=\Phi_0(x)+\int^t_0\hat b_s(Y_s)\dif s+\int^t_0\!\!\!\int_{|z|\leq 1}g_s(Y_{s-},z)\tilde N(\dif s,\dif z)
+\int^t_0\!\!\!\int_{|z|>1}g_s(Y_{s-},z) N(\dif s,\dif z),\label{Eq1}
\end{align}
where
\begin{align}
\tilde b_s(y):=\lambda \v_s(\Phi^{-1}_s(y))-\int_{|z|>1}[\v_s(\Phi^{-1}_s(y)+z)-\v_s(\Phi^{-1}_s(y))]\nu(\dif z)\label{Lp0}
\end{align}
and
\begin{align}
g_s(y,z):=\Phi_s(\Phi^{-1}_s(y)+z)-y.\label{Lp00}
\end{align}
\el
\begin{proof} Set
$$
\v^\eps_t(x):=(\v*\rho_\eps)(t,x),\ \ \Phi^\eps_t(x)=x+\v^\eps_t(x).
$$
By It\^o's formula, we have for all $t\in[0,T\wedge\tau)$,
\begin{align*}
\Phi^\eps_t(X_t)&=\Phi^\eps_0(X_0)+\int^t_0[\p_s\Phi^\eps_s(X_s)+( b^i_s\p_i\Phi^\eps_s)(X_s)]\dif s\\
&\quad+\int^t_0\!\!\!\int_{|z|\leq 1}[\Phi^\eps_s(X_{s-}+z)-\Phi^\eps_s(X_{s-})
-z^i\p_i\Phi^\eps_s(X_{s-})]\nu(\dif z)\dif s\\
&\quad+\int^t_0\!\!\!\int_{|z|\leq 1}[\Phi^\eps_s(X_{s-}+z)-\Phi^\eps_s(X_{s-})]\tilde N(\dif s,\dif z)\\
&\quad+\int^t_0\!\!\!\int_{|z|>1}[\Phi^\eps_s(X_{s-}+z)-\Phi^\eps_s(X_{s-})] N(\dif s,\dif z)\\
&=:\Phi^\eps_0(X_0)+I^\eps_1(t)+I^\eps_2(t)+I^\eps_3(t)+I^\eps_4(t).
\end{align*}
We want to take limits for the above equality. First of all,
for $I^\eps_4(t)$, by the dominated convergence theorem, we have
\begin{align*}
I^\eps_4(t)&=\int^t_0\!\!\!\int_{|z|>1}z N(\dif s,\dif z)+
\int^t_0\!\!\!\int_{|z|>1}[\v^\eps_s(X_{s-}+z)-\v^\eps_s(X_{s-})] N(\dif s,\dif z)\\
&\to\int^t_0\!\!\!\int_{|z|>1}z N(\dif s,\dif z)+
\int^t_0\!\!\!\int_{|z|>1}[\v_s(X_{s-}+z)-\v_s(X_{s-})] N(\dif s,\dif z)\\
&=\int^t_0\!\!\!\int_{|z|>1}[\Phi_s(X_{s-}+z)-\Phi_s(X_{s-})] N(\dif s,\dif z)
=\int^t_0\!\!\!\int_{|z|>1}g_s(Y_{s-},z) N(\dif s,\dif z),
\end{align*}
and for $I^\eps_3(t)$,
\begin{align*}
&\mE\left|\int^{t\wedge\tau}_0\!\!\!\int_{|z|\leq 1}[\Phi^\eps_s(X_{s-}+z)-\Phi^\eps_s(X_{s-})
-\Phi_s(X_{s-}+z)+\Phi_s(X_{s-})]\tilde N(\dif s,\dif z)\right|^2\\
&\quad=\mE\int^{t\wedge\tau}_0\!\!\!\int_{|z|\leq 1}|\Phi^\eps_s(X_{s-}+z)-\Phi^\eps_s(X_{s-})
-\Phi_s(X_{s-}+z)+\Phi_s(X_{s-})|^2\nu(\dif z)\dif s\to 0,
\end{align*}
where we have used that for some $C$ independent of $\eps$,
$$
|\Phi^\eps_s(X_{s-}+z)-\Phi^\eps_s(X_{s-})|\leq C|z|^2.
$$
Noting that
$$
\p_s\Phi^\eps=\p_s\v^\eps=-(\cL_0-\lambda) \v^\eps-(b^i\p_i \v)*\rho_\eps-b*\rho_\eps
=-(\cL_0-\lambda) \v^\eps-(b^i\p_i \Phi)*\rho_\eps,
$$
we have
\begin{align*}
I^\eps_1(t)+I^\eps_2(t)&=\lambda\int^t_0\v^\eps_s(X_s)\dif s-\int^t_0\!\!\!\int_{|z|>1}
[\v^\eps_s(X_s+z)-\v^\eps_s(X_s)]\nu(\dif z)\dif s\\
&\quad+\int^t_0 [(b^i_s\p_i\Phi^\eps_s)(X_s)-((b^i\p_i \Phi)*\rho_\eps)(s,X_s)]\dif s.
\end{align*}
By the dominated convergence theorem, the first two terms converge to
$$
\lambda\int^t_0\v_s(X_s)\dif s-\int^t_0\!\!\!\int_{|z|>1}
[\v_s(X_s+z)-\v_s(X_s)]\nu(\dif z)\dif s=\int^t_0\tilde b_s(Y_s)\dif s.
$$
Using Krylov's estimate (\ref{Lp66}), we have
\begin{align*}
&\mE\int^{t\wedge\tau}_0|(b^i_s\p_i\Phi^\eps_s)(X_s)-((b^i\p_i \Phi)*\rho_\eps)(s,X_s)|\dif s\\
&\qquad\leq C\int^t_0\left(\int_{\mR^d}|(b^i_s\p_i\Phi^\eps_s)(x)-((b^i\p_i \Phi)*\rho_\eps)(s,x)|^p\dif x\right)^{\frac{q}{p}}
\dif s\to 0,
\end{align*}
where $q>\frac{\alpha}{\alpha-1}$. Combining the above calculations, we obtain that $Y_t$ solves (\ref{Eq1}).
\end{proof}

We are now in a position to give:

\begin{proof}[Proof of Theorem \ref{Main}]
We first assume that for some $\beta,p$ satisfying (\ref{Con}),
$$
b\in L^\infty_{loc}(\mR_+; L^\infty(\mR^d)\cap\mW^\beta_p).
$$
The existence of weak solutions has been obtained in Theorem \ref{Weak}.
Below, we concentrate on the proof of the pathwise uniqueness.

Fix an ($\sF_t$)-stopping time $\tau$ and let
$X_t, \hat X_t\in\sS^{\tau}_b(x)$ be two solutions of SDE (\ref{SDE}).
Fixing $T>0$, we want to prove that
$$
Y_t:=\Phi_t(X_t)=\Phi_t(\hat X_t)=:\hat Y_t,\ \ \forall t\in[0,T\wedge\tau).
$$
Define $\sigma_0\equiv 0$ and for $n\in\mN$,
$$
\sigma_n:=\inf\{t\geq \sigma_{n-1}: |L_t-L_{t-}|>1\}.
$$
Set
$$
\sigma^T_n=\sigma_n\wedge T\wedge\tau.
$$
Recall (\ref{Lev}) and
$$
\int^t_0\!\!\!\int_{|z|>1}g_s(Y_{s-},z) N(\dif s,\dif z)=\sum_{s\in(0,t]}g_s(Y_{s-},L_s-L_{s-})\cdot 1_{|L_s-L_{s-}|>1}.
$$
By Lemma \ref{Le3}, $Z_t:=Y_t-\hat Y_t$ satisfy the following equation on random interval
$[\sigma^T_{n},\sigma^T_{n+1})$:
\begin{align}
Z_t=Z_{\sigma^T_{n}}+\int^t_{\sigma^T_{n}}[\tilde b_s(Y_s)-\tilde b_s(\hat Y_s)]\dif s
+\int^t_{\sigma^T_{n}}\!\int_{|z|\leq1}[g_s(Y_{s-},z)-g_s(\hat Y_{s-},z)]\tilde N(\dif s,\dif z).\label{Op9}
\end{align}
Let us first prove that
$$
Z_t=0\mbox{ a.s. on $[0,\sigma^T_1)$.}
$$
Note that by  (\ref{Lp0}), (\ref{Lp8})  and (\ref{Op1}),
\begin{align}
|\tilde b_s(y)-\tilde b_s(y')|\leq C|y-y'|,\label{Op2}
\end{align}
and by (\ref{Es2}),
\begin{align}
|g_s(y,z)-g_s(y',z)|&=|\Phi_s(\Phi^{-1}_s(y)+z)-\Phi_s(\Phi^{-1}_s(y))-\Phi_s(\Phi^{-1}_s(y')+z)+\Phi_s(\Phi^{-1}_s(y')|\no\\
&=|\v_s(\Phi^{-1}_s(y)+z)-\v_s(\Phi^{-1}_s(y))-\v_s(\Phi^{-1}_s(y')+z)+\v_s(\Phi^{-1}_s(y')|\no\\
&=|(\sT_z\v_s)(\Phi^{-1}_s(y))-(\sT_z\v_s)(\Phi^{-1}_s(y'))|\no\\
&\leq C|\Phi^{-1}_s(y)-\Phi^{-1}_s(y')|\cdot(\cM|\nabla\sT_z\v_s|(\Phi^{-1}_s(y))+
\cM|\nabla\sT_z\v_s|(\Phi^{-1}_s(y')))\no\\
&\leq C|y-y'|\cdot(\cM|\nabla\sT_z\v_s|(\Phi^{-1}_s(y))+
\cM|\nabla\sT_z\v_s|(\Phi^{-1}_s(y'))).\label{Op3}
\end{align}
Since $\mE|X_t|^2=+\infty$, for taking expectations for (\ref{Op9}),
we need to use stopping time to cut off it. For $R>0$, define
\begin{align}
\zeta_R:=\inf\{t\geq0: |X_t|\vee|\hat X_t|\geq R\}.\label{Cut}
\end{align}
Let $\eta$ be any ($\sF_t$)-stopping time. By (\ref{Op9}), (\ref{Op2}) and (\ref{Op3}), we have
\begin{align*}
\mE|Z_{t\wedge\sigma^T_1\wedge\zeta_R\wedge\eta-}|^2&\leq
C\mE\int^{t\wedge\sigma^T_1\wedge\zeta_R\wedge\eta}_0\left(|Z_s|^2
+\int_{|z|\leq1}|g_s(Y_{s-},z)-g_s(\hat Y_{s-},z)|^2\nu(\dif z)\right)\dif s\\
&\leq C\mE\int^{t\wedge\sigma^T_1\wedge\zeta_R\wedge\eta}_0|Z_{s-}|^2\dif (s+A_s)
\leq C\mE\int^{t\wedge\eta}_0|Z_{s\wedge\sigma^T_1\wedge\zeta_R-}|^2\dif (s+A_{s\wedge\tau}),
\end{align*}
where
$$
A_t:=\int^t_0\!\!\!\int_{|z|\leq 1}\Big(\cM|\nabla\sT_z\v_s|(\Phi^{-1}_s(Y_s))
+\cM|\nabla\sT_z\v_s|(\Phi^{-1}_s(\hat Y_s))\Big)^2\nu(\dif z)\dif s.
$$
By Fubini's theorem, we have
\begin{align*}
\mE A_{t\wedge\tau}&=\int_{|z|\leq 1}\mE\int^{t\wedge\tau}_0\Big(\cM|\nabla\sT_z\v_s|(X_s)
+\cM|\nabla\sT_z\v_s|(\hat X_s)\Big)^2\dif s\nu(\dif z)\\
&\stackrel{(\ref{Lp66})}{\leq} C\int_{|z|\leq 1}\sup_{s\in[0,t]}\|(\cM|\nabla\sT_z\v_s|)^2\|_{p/2}\nu(\dif z)\\
&=C\int_{|z|\leq 1}\sup_{s\in[0,t]}\|\cM|\nabla\sT_z\v_s|\|^2_p\nu(\dif z)
\stackrel{(\ref{Es30})}{\leq} C\int_{|z|\leq 1}\sup_{s\in[0,t]}\|\sT_z\v_s\|^2_{1,p}\nu(\dif z)\\
&\stackrel{(\ref{Lp7})}{\leq} C\sup_{s\in[0,t]}\|\v_s\|^2_{\gamma+\beta,p}
\int_{|z|\leq 1}|z|^{2(\gamma+\beta-1)}\nu(\dif z)<+\infty,
\end{align*}
where in the last step we have used (\ref{Lp8}), $2(\gamma+\beta-1)>\alpha$ and  (\ref{Coa}).
Therefore, $t\mapsto A_{t\wedge\tau}$ is a continuous ($\sF_t$)-adapted increasing process.
By Lemma \ref{Le0}, we obtain that for all $t\geq 0$,
$$
Z_{t\wedge\sigma^T_1\wedge\zeta_R-}=0,\ \ a.s.
$$
Letting $R\to\infty$ yields that  for all $t\in[0,T\wedge\tau)$,
$$
Z_{t\wedge\sigma_1-}=Z_{t\wedge\sigma^T_1-}=0,\ \ a.s.
$$
Thus, if $\sigma_1<T\wedge\tau$, then
$$
Z_{\sigma_1}=Z_{\sigma_1-}+[g_{\sigma_1}(Y_{\sigma_1-},L_{\sigma_1}-L_{\sigma_1-})
-g_{\sigma_1}(\hat Y_{\sigma_1-},L_{\sigma_1}-L_{\sigma_1-})]=0.
$$
Repeating the above calculations, we find that for all $n\in\mN$ and $t\in[0,T\wedge\tau)$,
$$
Z_{t\wedge\sigma_n-}=0\ \ a.s.
$$
Letting $n, T\to\infty$ produces that for all $t\in[0,\tau)$,
$$
Z_t=0\Rightarrow Y_t=\hat Y_t\Rightarrow X_t=\hat X_t\ \ a.s.
$$

Lastly, we assume that $b$ satisfies (\ref{Op6}) and (\ref{Op06}). For $n\in\mN$, let $\chi_n\in C^\infty_0(\mR^d)$ with
$\chi_n(x)=1$ for $|x|\leq n$ and $\chi_n(x)=0$ for $|x|>n+1$. Define
$$
b^{(n)}_t(x)=b_t(x)\chi_n(x).
$$
Then $b^{(n)}\in L^\infty_{loc}(\mR_+; L^\infty(\mR^d)\times\mW^\beta_p)$.
By the previous proof, for each $x\in\mR^d$,
there exists a unique strong solution $X^{(n)}_t\in\sS^\infty_{b^{(n)}}(x)$
to SDE (\ref{SDE}) with drift $b^{(n)}$. For $n\geq k$, define
$$
\tau_{n,k}(x,\omega):=\inf\{t\geq 0: |X^{(n)}_t(\omega,x)|\geq k\}.
$$
It is easy to see that
$$
X^{(n)}_t(x), X^{(k)}_t(x)\in\sS^{\tau_{n,k}(x)}_{b^{(k)}}(x).
$$
As the local uniqueness has been proven, we have
$$
P\{\omega:X^{(n)}_t(\omega,x)= X^{(k)}_t(\omega,x), \forall t\in[0,\tau_{n,k}(x,\omega))\}=1,
$$
which implies that for $n\geq k$,
$$
\tau_{k,k}(x)\leq\tau_{n,k}(x)\leq\tau_{n,n}(x), \ \ a.s.
$$
Hence, if we let $\zeta_k(x):=\tau_{k,k}(x)$, then $\zeta_k(x)$ is an increasing sequence of
($\sF_t$)-stopping times and for $n\geq k$,
$$
P\{\omega: X^{(n)}_t(x,\omega)=X^{(k)}_t(x,\omega),\ \ \forall t\in[0,\zeta_k(x,\omega))\}=1.
$$
Now, for each $k\in\mN$, we can define $X_t(x,\omega)=X^{(k)}_t(x,\omega)$ for
$t<\zeta_k(x,\omega)$ and $\zeta(x)=\lim_{k\to\infty}\zeta_k(x)$.
It is clear that $X_t(x)\in\sS^{\zeta(x)}_{b}(x)$ and (\ref{HH1}) holds.
\end{proof}

{\bf Acknowledgements:}

The author would like to express his deep thanks to Professor Jiagang Ren for teaching him the interpolation techniques
and the fact (\ref{Fact}) about ten years ago. The supports of NSFs of China (No. 10971076; 10871215) are also
acknowledged.


\begin{thebibliography}{999}

\bibitem{Ad}Aldous, D.: Stopping times and tightness. Ann. Prob., Vol. 6(1978), 335-340.

\bibitem{Ar-Zh}Airault, H. and Zhang, X.: Smoothness of indicator
functions of some sets in Wiener spaces.  J. Math. Pures Appl., Vol. 79, No. 5, (2000) 515-523.

\bibitem{Ap}Applebaum D. L\'evy processes and stochastic calculus. Cambridge Studies in Advance Mathematics 93, Cambridge University PRess, 2004.

\bibitem{Ba1}Bass, R.: Stochastic differential equations driven by symmetric stable processes.
Seminaire de Probabilities, XXXVI, Lecture Notes in Math., 1801,302-313(2003).

\bibitem{Ba3} Bass, R.: Stochastic differential equations with jumps.
Probability Surveys, Vol.1, (2004)1-19.

\bibitem{Bo-Ja}Bogdan K. and Jakubowski T.: Estimates of heat kernel of fractional Laplacian perturbed by
gradient operators. Commun. Math. Phys., 271(2007), 179-198.

\bibitem{Ch-Ki-So}Chen Z., Kim P. and Song R.: Dirichlet heat kernel estimates for fractional Laplacian with gradient
perturbation. http://arxiv.org/abs/1011.3273.

\bibitem{Pr}Priola E.: Pathwise uniqueness for singular SDEs driven by stable processes.
http://arxiv.org/abs/1005.4237


\bibitem{Ik-Wa}Ikeda N., Watanabe, S.: Stochastic differential equations
and diffusion processes, 2nd ed., North-Holland/Kodanska, Amsterdam/Tokyo, 1989.

\bibitem{Kr}Krylov N.V.: Controlled diffusion processes. Translated from
the Russian by A.B. Aries. Applications of Mathematics, 14. Springer-Verlag,
New York-Berlin, 1980.


\bibitem{Kr-Ro}Krylov N.V. and R\"ockner M.: Strong solutions of stochasitc equations
with singluar time dependent drift. Probab. Theory Relat. Fields, 131, 154-196(2005).

\bibitem{Ku}Kurenok V.P.: Stochastic equations with time time-dependent drifts driven by L\'evy processes.
J. Theortic Prob., 20(2007)859-869.

\bibitem{Ku1}Kurenok V.P.: A note on $L_2$-estimates for stable integrals with drift.
Trans. Amer. Math. Soc. 360 (2008), 925-938.


\bibitem{Re-Zh}Ren, J. and Zhang, X.: Limit theorems for stochastic differential equations
with discontinuous coefficients. to appear in SIAM J. Math. Anal.

\bibitem{Sa}Sato K.I.: L\'evy processes and infinite divisble distributions. Cambridge University Press, Cambridge, 1999.

\bibitem{St}Stein E.M.:  Singular integrals and differentiability properties of functions.
Princeton, N.J.,  Princeton University Press,  1970.

\bibitem{Ta-Ts-Wa}Tanaka H., Tsuchiya M. and Watanabe S.: Perturbation of drift-type for L\'evy processes. J. Math. Kyoto Univ. 14(1974), 73-92.

\bibitem{Tr}Triebel, H.: Interpolation Theory, Function Spaces, Differential Operators. North-Holland Publishing Company,
Amsterdam, 1978.

\bibitem{Ve}Veretennikov, A. Ju.: On the strong solutions of stochastic differential
equations. Theory Probab. Appl., 24(1979),354-366.

\bibitem{Zh0}Zhang X.: Stochastic homeomorphism flows of SDEs with singular drifts and Sobolev diffusion coefficients.
http://arxiv.org/abs/1010.3403.

\bibitem{Zh1}Zhang X.: Well-posedness and large deviation for degenerate SDEs with Sobolev coefficients.
http://arxiv.org/abs/1002.4297.

\bibitem{Zv}Zvonkin, A.K.: A transformation of the phase space of a diffusion process
that removes the drift.  Mat. Sbornik, No.1, 93(135),129-149(1974).

\end{thebibliography}
\end{document}